\documentclass[graybox]{svmult}
\usepackage{newtxtext,newtxmath}

\newcommand{\wto}{\rightharpoonup}
\newcommand{\cptemb}{\hookto\hspace*{-0.85em}\hookto}
\newcommand{\setb}[2]{\big\{#1\,:\,#2\big\}}
\newcommand{\ovl}[1]{\overline{#1}}
\newcommand{\hookto}{\hookrightarrow}
\newcommand{\reals}{\mathbb{R}}
\DeclareMathOperator{\supp}{supp}
\DeclareMathOperator{\dist}{dist}
\newcommand{\rt}{\reals^{3}}
\newcommand{\om}{\Omega}
\newcommand{\ga}{\Gamma}
\newcommand{\gat}{\ga_{\mathsf{t}}}
\newcommand{\gan}{\ga_{\mathsf{n}}}
\DeclareMathOperator{\p}{\partial}
\DeclareMathOperator{\na}{\nabla}
\DeclareMathOperator{\rot}{rot}
\DeclareMathOperator{\curl}{curl}
\DeclareMathOperator{\divergence}{div}
\renewcommand{\div}{\divergence}
\newcommand{\lsymbol}{\mathsf{L}}
\newcommand{\lgen}[3]{\overset{#1}{\lsymbol}{}^{#2}_{#3}}
\newcommand{\lt}{\lgen{}{2}{}}
\newcommand{\ltom}{\lt(\om)}
\newcommand{\csymbol}{\mathsf{C}}
\newcommand{\cgen}[3]{\overset{#1}{\csymbol}{}^{#2}_{#3}}
\newcommand{\ci}{\cgen{}{\infty}{}}
\newcommand{\hsymbol}{\mathsf{H}}
\newcommand{\hgen}[3]{\overset{#1}{\hsymbol}{}^{#2}_{#3}}
\newcommand{\hoom}{\ho(\om)}
\newcommand{\ho}{\hgen{}{1}{}}
\newcommand{\hoc}{\mathring\hsymbol^1}
\newcommand{\hocgat}{\hoc_{\gat}}
\newcommand{\hocgatom}{\hocgat(\om)}
\newcommand{\hmo}{\hsymbol^{-1}}
\newcommand{\hmoom}{\hmo(\om)}
\newcommand{\hmoc}{\mathring\hsymbol^{-1}}
\newcommand{\hmocom}{\hmoc(\om)}
\newcommand{\rsymbol}{\mathsf{R}}
\newcommand{\rgen}[3]{\overset{#1}{\rsymbol}{}^{#2}_{#3}}
\newcommand{\rc}{\mathring\rsymbol}
\newcommand{\rcgat}{\rc_{\gat}}
\newcommand{\rcgatom}{\rcgat(\om)}
\newcommand{\dsymbol}{\mathsf{D}}
\newcommand{\dgen}[3]{\overset{#1}{\dsymbol}{}^{#2}_{#3}}
\newcommand{\dcgatom}{\dcgat(\om)}
\newcommand{\dc}{\mathring\dsymbol}
\newcommand{\dcgat}{\dc_{\gat}}
\renewcommand{\H}{\mathsf{H}}
\newcommand{\harmsymbol}{\mathcal{H}}
\newcommand{\harmgen}[3]{\overset{#1}{\harmsymbol}{}^{#2}_{#3}}
\newcommand{\harm}{\harmgen{}{}{}}
\newcommand{\harmom}{\harm(\om)}
\newcommand{\cm}{c_{\mathsf{m}}}
\newcommand{\cfp}{c_{\mathsf{f,p}}}
\newcommand{\norm}[1]{|#1|}
\newcommand{\normltom}[1]{\norm{#1}_{\ltom}}
\newcommand{\scp}[2]{\langle#1,#2\rangle}
\newcommand{\scpltom}[2]{\scp{#1}{#2}_{\ltom}}

\begin{document}

\title*{A Global div-curl-Lemma for Mixed Boundary Conditions in Weak Lipschitz Domains}
\titlerunning{A Global div-curl-Lemma for Mixed Boundary Conditions in Weak Lipschitz Domains} 
\author{Dirk Pauly}
\authorrunning{Dirk Pauly}
\institute{Dirk Pauly \at 
Fakult\"at f\"ur Mathematik, Universit\"at Duisburg-Essen, Campus Essen, Germany\\
\email{dirk.pauly@uni-due.de}}

\maketitle

\vspace*{-28mm}
\abstract{We prove a global version of the so-called $\div$-$\curl$-lemma,
a crucial result for compensated compactness and in homogenization theory,
for mixed tangential and normal boundary conditions in bounded weak Lipschitz domains in 3D
and weak Lipschitz interfaces.
The crucial tools and the core of our arguments are 
the de Rham complex and Weck's selection theorem,
the essential compact embedding result for Maxwell's equations.}

\section{Introduction and Main Results}
\label{introsec}

We shall prove a global (and hence also a local) version of the so-called $\div$-$\curl$-lemma,
with mixed tangential and normal boundary conditions 
for bounded weak Lipschitz domains $\om$ in 3D, more precisely
for admissible pairs $(\om,\gat)$
of a bounded weak Lipschitz domain $\om\subset\rt$ 
and a part $\gat$ of its boundary $\ga$, see Definition \ref{div-rot-lem-def} for details.

\begin{theorem}[global $\div$-$\curl$-lemma]
Let $(\om,\gat)$ be admissible and let
\begin{itemize}
\item[\bf(i)]\quad
$E_{n},E\in D(\curl_{\gat})$, 
\item[\bf(i')]\quad
$E_{n}\wto E$ in\footnote{In particular, 
$E_{n}\wto E$ in $\ltom$ and $\curl E_{n}\wto \curl E$ in $\ltom$.}
$D(\curl_{\gat})$,
\item[\bf(ii)]\quad
$H_{n},H\in D(\div_{\gan})$, 
\item[\bf(ii')]\quad
$H_{n}\wto H$ in\footnote{In particular, 
$H_{n}\wto H$ in $\ltom$ and $\div H_{n}\wto \div H$ in $\ltom$.}
$D(\div_{\gan})$.
\end{itemize}
Then
\begin{itemize}
\item[\bf(iii)]\quad
$\scpltom{E_{n}}{H_{n}}\to\scpltom{E}{H}$.
\end{itemize}
\label{div-rot-lem}
\end{theorem}

Here, we introduce the densely defined and closed linear operators
$\na_{\gat}$, $\curl_{\gat}$, $\div_{\gat}$
with domains of definition\footnote{Various notations like
\begin{align*}
D(\na_{\gat})&=\H(\na_{\gat},\om)=\H_{\gat}(\na,\om)=\H^{1}_{\gat}(\om)=\hocgatom,\\
D(\curl_{\gat})&=\H(\curl_{\gat},\om)=\H_{\gat}(\curl,\om)=\rgen{}{}{\gat}(\om)=\rcgatom,\\
D(\div_{\gat})&=\H(\div_{\gat},\om)=\H_{\gat}(\div,\om)=\dgen{}{}{\gat}(\om)=\dcgatom
\end{align*}
can be found frequently in the literature, where also $\curl=\rot$ is used.} 
$D(\na_{\gat})$, $D(\curl_{\gat})$, $D(\div_{\gat})$
as closures of the classical differential operators from vector analysis acting on $\ltom$ and
defined on smooth test functions resp. test vector fields
bounded away from the boundary part $\gat$ given by
$$\cgen{}{\infty}{\gat}(\om):=
\setb{\varphi|_\om}
{\varphi\in\ci(\rt),\,\supp\varphi\text{ compact},\,\dist(\supp\varphi,\gat)>0}.$$
As shown in \cite[Theorem 4.5]{bauerpaulyschomburgmcpweaklip} 
(weak equals strong in terms of definitions of boundary conditions)
their adjoints are given by $-\div_{\gan}$, $\curl_{\gan}$, $-\na_{\gan}$
defined in the same way. Note that these operators are unbounded
and that the domains of definition are Hilbert spaces equipped with the 
respective proper graph inner products.

\begin{corollary}[local $\div$-$\curl$-lemma]
Let $\om\subset\rt$ be an open set and let
\begin{itemize}
\item[\bf(i)]\quad
$E_{n},E\in D(\curl)$,
\item[\bf(i')]\quad
$E_{n}\wto E$ in $D(\curl)$,
\item[\bf(ii)]\quad
$H_{n},H\in D(\div)$, 
\item[\bf(ii')]\quad
$H_{n}\wto H$ in $D(\div)$.
\end{itemize}
Then
\begin{itemize}
\item[\bf(iii)]\quad
$\forall\,\varphi\in\cgen{}{\infty}{\ga}(\om)\qquad\scpltom{\varphi\,E_{n}}{H_{n}}\to\scpltom{\varphi\,E}{H}$.
\end{itemize}
$\cgen{}{\infty}{\ga}(\om)$ 
may be replaced by $\cgen{}{1}{\ga}(\om)$ or even $\cgen{}{0,1}{\ga}(\om)$,
the space of Lipschitz continuous functions vanishing in a neighbourhood of $\ga$.
Moreover, the boundedness of $(E_{n})$ and $(H_{n})$ in local spaces is sufficient 
for the assertion to hold.
\label{div-rot-lem-loc}
\end{corollary}

The $\div$-$\curl$-lemma, or compensated compactness,
see the original papers by Murat \cite{murat1978} and Tartar \cite{tartar1979}
or \cite{evans1990,struwe2008},
and its variants and extensions have plenty important applications.
It is widely used in the theory of homogenization of (nonlinear) partial differential equations,
see, e.g., \cite{struwe2008}.
Moreover, it is crucial in establishing compactness and regularity results 
for nonlinear partial differential equations such as harmonic maps,
see, e.g., \cite{freiremuellerstruwe1998,evans1991,riviere2007}.
Numerical applications can be found, e.g., in \cite{bartels2010}.
The $\div$-$\curl$-lemma is further a crucial tool 
in the homogenization of stochastic partial differential equations,
especially with certain random coefficients, see, e.g.,
the survey \cite{alexanderian2015}
and the literature cited therein, e.g., \cite{glorianeukammotto2015}. 

For an extensive discussion and a historical overview of the $\div$-$\curl$-lemma see \cite{tartar2009}.
More recent discussions can be found, e.g., in \cite{brianecasadodiazmurat2009,tartar2015}.
Recently, in \cite{waurick2018a,paulydivcurl}
the $\div$-$\curl$-lemma has been proved in a general Hilbert space setting
which allows for various applications in mathematical physics.
Interesting and new applications to homogenization of partial differential equations 
can be found in \cite{waurick2018b}.

Let us also mention that the $\div$-$\curl$-lemma
is particularly useful to treat homogenization of problems arising in plasticity,
see, e.g., a recent preprint on this topic \cite{roegerschweizer2017},
for which the preprint \cite{schweizer2017} provides the important key $\div$-$\curl$-lemma.
Unfortunately, in \cite{schweizer2017,roegerschweizer2017} 
a $\hoom$-detour is used as the core argument for the proofs.
The same detour is utilized in the recent contribution \cite{kozonoyanagisawa2013} 
where $\div$-$\curl$-type lemmas are presented
which also allow for inhomogeneous boundary conditions.
This unnecessarily high regularity assumption of $\hoom$-fields
excludes results like
\cite{kozonoyanagisawa2013,schweizer2017,roegerschweizer2017}
to be applied to important applications 
which are stated, e.g., in Lipschitz domains.

\section{Notations, Preliminaries, and Proofs}
\label{prelimsec}

\begin{definition}[admissible domains]
We call a pair $(\om,\gat)$ admissible, if 
\begin{itemize}
\item[\bf(i)]\quad
$\om\subset\rt$ is a bounded weak Lipschitz domain 
in the sense of \cite[Definition 2.3]{bauerpaulyschomburgmcpweaklip}
\item[\bf(ii)]\quad
with boundary $\ga:=\p\om$, which is divided into two relatively open weak Lipschitz subsets 
$\gat\subset\ga$ and its complement $\gan:=\ga\setminus\ovl{\gat}$
in the sense of \cite[Definition 2.5]{bauerpaulyschomburgmcpweaklip}.
\end{itemize}
\label{div-rot-lem-def}
\end{definition}

Note that strong Lipschitz domains (locally below a graph of a Lipschitz function) 
are weak Lipschitz domains (the boundary is a Lipschitz manifold)
which holds for the boundary as well as for the interface. 
The reverse implication is not true due to the failure of the
implicit function theorem for Lipschitz mappings.
Throughout this paper we shall assume the latter regularity of $\om$, and $\ga$, $\gat$, $\gan$.

Recently, in \cite{bauerpaulyschomburgmcpweaklip}, 
Weck's selection theorem \cite{weckmax}, also known as the Maxwell compactness property, has been shown
to hold for such bounded weak Lipschitz domains and mixed boundary conditions.
More precisely, the following holds:

\begin{lemma}[Weck's selection theorem]
Let $(\om,\gat)$ be admissible. Then the embedding
$$D(\curl_{\gat})\cap D(\div_{\gan})\cptemb\ltom$$
is compact.
\label{weckstlem}
\end{lemma}

For a proof see \cite[Theorem 4.7]{bauerpaulyschomburgmcpweaklip}. 
A short historical overview of Weck's selection theorem
is given in the introduction of \cite{bauerpaulyschomburgmcpweaklip},
see also the original paper \cite{weckmax} and 
\cite{picardcomimb,webercompmax,costabelremmaxlip,witschremmax,jochmanncompembmaxmixbc,leisbook}
for simpler proofs and generalizations.

Let us emphasize that our assumptions also allow for Rellich's selection theorem, i.e., the embedding
\begin{align}
\label{rellichst}
D(\na_{\gat})\cptemb\ltom
\end{align}
is compact, see, e.g., \cite[Theorem 4.8]{bauerpaulyschomburgmcpweaklip}. 
By density we have the two rules of integration by parts
\begin{align}
\label{partintna}
\forall\,u&\in D(\na_{\gat})
&
\forall\,H&\in D(\div_{\gan})
&
\scpltom{\na u}{H}&=-\scpltom{u}{\div H},\\
\label{partintrot}
\forall\,E&\in D(\curl_{\gat})
&
\forall\,H&\in D(\curl_{\gan})
&
\scpltom{\curl E}{H}&=\scpltom{E}{\curl H}.
\end{align}

A direct consequence of Lemma \ref{weckstlem} is the compactness of the unit ball in 
$$\harmom:=N(\curl_{\gat})\cap N(\div_{\gan}),$$
the space of so-called Dirichlet-Neumann fields. Hence $\harmom$ is finite-dimensional.
Here and in the following we denote the kernels and the ranges of our operators
$\na_{\gat}$, $\curl_{\gat}$, $\div_{\gat}$ by
$$N(\na_{\gat}),\quad
N(\curl_{\gat}),\quad
N(\div_{\gat}),\qquad
R(\na_{\gat}),\quad
R(\curl_{\gat}),\quad
R(\div_{\gat}).$$
Another immediate consequence of Weck's selection theorem, Lemma \ref{weckstlem}, 
using a standard indirect argument, is the so-called Maxwell estimate, i.e.,
there exists $\cm>0$ such that for all 
$E\in D(\curl_{\gat})\cap D(\div_{\gan})\cap\harmom^{\bot_{\ltom}}$
\begin{align}
\label{maxest}
\normltom{E}\leq\cm\big(\normltom{\curl E}+\normltom{\div E}\big),
\end{align}
see \cite[Theorem 5.1]{bauerpaulyschomburgmcpweaklip}.
Recent estimates for the Maxwell constant $\cm$ can be found in \cite{paulymaxconst0,paulymaxconst1,paulymaxconst2}.
Analogously, Rellich's selection theorem \eqref{rellichst} shows the Friedrichs/Poincar\'e estimate, i.e.,
there exists $\cfp>0$ such that for all $u\in D(\na_{\gat})$
\begin{align}
\label{fpest}
\normltom{u}\leq\cfp\normltom{\na u},
\end{align}
see \cite[Theorem 4.8]{bauerpaulyschomburgmcpweaklip}.
To avoid case studies due to the one-dimensional kernel $\reals$ of $\na$
when using the Friedrichs/Poincar\'e estimate in the case $\gat=\emptyset$, we also define 
$$D(\na_{\emptyset})
:=D(\na)\cap\reals^{\bot_{\ltom}}
=\setb{u\in\hoom}{\int_{\om}u=0}.$$
By the projection theorem, applied to our densely defined and closed (unbounded) linear operator
$$\na_{\gat}:D(\na_{\gat})\subset\ltom\longrightarrow\ltom$$
with (Hilbert space) adjoint
$$\na_{\gat}^{*}=-\div_{\gan}:D(\div_{\gan})\subset\ltom\longrightarrow\ltom,$$
where we have used \cite[Theorem 4.5]{bauerpaulyschomburgmcpweaklip} (weak equals strong),
we get the simple (orthogonal) Helmholtz decomposition
\begin{align}
\label{helmsim}
\ltom
&=R(\na_{\gat})\oplus_{\ltom}N(\div_{\gan}),
\end{align}
see \cite[Theorem 5.3 or (13)]{bauerpaulyschomburgmcpweaklip},
which immediately implies the orthogonal decomposition
\begin{align}
\label{helmsimtwo}
D(\curl_{\gat})
&=R(\na_{\gat})\oplus_{\ltom}\big(D(\curl_{\gat})\cap N(\div_{\gan})\big)
\end{align}
as the complex property $R(\na_{\gat})\subset N(\curl_{\gat})$ holds. 
Here $\oplus_{\ltom}$ in the decompositions \eqref{helmsim} and \eqref{helmsimtwo}
denotes the orthogonal sum in the Hilbert space $\ltom$.
By \eqref{fpest} the range $R(\na_{\gat})$ is closed in $\ltom$,
see also \cite[Lemma 5.2]{bauerpaulyschomburgmcpweaklip}.
Note that we call \eqref{helmsim} a simple Helmholtz decomposition, 
since the refined Helmholtz decomposition
$$\ltom=R(\na_{\gat})\oplus_{\ltom}\harmom\oplus_{\ltom}R(\curl_{\gan})$$
holds as well, see \cite[Theorem 5.3]{bauerpaulyschomburgmcpweaklip},
where also $R(\curl_{\gan})$ is closed in $\ltom$
as a consequence of \eqref{maxest}, 
see \cite[Lemma 5.2]{bauerpaulyschomburgmcpweaklip}.

\vspace*{3mm}\noindent
{\bf Proof of Theorem \ref{div-rot-lem}}
By \eqref{helmsimtwo} we have $D(\curl_{\gat})\ni E_{n}=\na u_{n}+\tilde{E}_{n}$
with some $u_{n}\in D(\na_{\gat})$
and $\tilde{E}_{n}\in D(\curl_{\gat})\cap N(\div_{\gan})$. 
Then $(u_{n})$ is bounded in $\hoom$ by orthogonality 
and the Friedrichs/Poincar\'e estimate \eqref{fpest}.
By orthogonality $(\tilde{E}_{n})$ is bounded in $D(\curl_{\gat})\cap N(\div_{\gan})$
and $\curl\tilde{E}_{n}=\curl E_{n}$.
Hence, using Rellich's and Weck's selection theorems 
there exist $u\in D(\na_{\gat})$ and $\tilde{E}\in D(\curl_{\gat})\cap N(\div_{\gan})$
and we can extract two subsequences, again denoted by $(u_{n})$ and $(\tilde{E}_{n})$
such that $u_{n}\rightharpoonup u$ in $D(\na_{\gat})$ and $u_{n}\to u$ in $\ltom$
as well as $\tilde{E}_{n}\rightharpoonup\tilde{E}$ in $D(\curl_{\gat})\cap N(\div_{\gan})$ 
and $\tilde{E}_{n}\to\tilde{E}$ in $\ltom$. 
We observe $E=\na u+\tilde{E}$, giving the simple Helmholtz decomposition for $E$.
Finally, by \eqref{partintna}
\begin{align*}
\scpltom{E_{n}}{H_{n}}
&=\scpltom{\na u_{n}}{H_{n}}
+\scpltom{\tilde{E}_{n}}{H_{n}}\\
&=-\scpltom{u_{n}}{\div H_{n}}
+\scpltom{\tilde{E}_{n}}{H_{n}}\\
&\to-\scpltom{u}{\div H}
+\scpltom{\tilde{E}}{H}\\
&=\scpltom{\na u}{H}
+\scpltom{\tilde{E}}{H}\\
&=\scpltom{E}{H}.
\end{align*}
As the limit is unique, the original sequence $\big(\scpltom{E_{n}}{H_{n}}\big)$ already converges
to the limit $\scpltom{E}{H}$.
$\blacksquare$

\vspace*{3mm}\noindent
{\bf Proof of Corollary \ref{div-rot-lem-loc}}
Let $\gat:=\ga$ and hence $\gan=\emptyset$.
$(\varphi\,E_{n})$ is bounded in $D(\curl_{\ga})$
and $(H_{n})$ is bounded in $D(\div)$.
Theorem \ref{div-rot-lem} shows the assertion.
$\blacksquare$

\section{Generalizations and the Classical div-curl-Lemma}
\label{gensec}

In \cite{waurick2018a,paulydivcurl} more general $\div$-$\curl$-lemmas have been presented.
In particular in \cite{paulydivcurl} we can find the following generalization
to distributions.

\begin{theorem}[alternative global $\div$-$\curl$-lemma]
Let $\om\subset\rt$ be a bounded strong Lipschitz domain with trivial topology and let
\begin{itemize}
\item[\bf(i)]\quad
$E_{n},H_{n},E,H\in\ltom$, 
\item[\bf(i')]\quad
$E_{n}\wto E$ and $H_{n}\wto H$ in $\ltom$.
\end{itemize}
Moreover, let either
\begin{itemize}
\item[\bf(ii)]\quad
$(\widehat{\curl}E_{n})$ be relatively compact in $\hmocom$,
\item[\bf(iii)]\quad
$(\widetilde{\div}H_{n})$ be relatively compact in $\hmoom$,
\end{itemize}
or
\begin{itemize}
\item[\bf(ii')]\quad
$(\widetilde{\curl}E_{n})$ be relatively compact in $\hmoom$,
\item[\bf(iii')]\quad
$(\widehat{\div}H_{n})$ be relatively compact in $\hmocom$.
\end{itemize}
Then
\begin{itemize}
\item[\bf(iv)]\quad
$\scpltom{E_{n}}{H_{n}}\to\scpltom{E}{H}$.
\end{itemize}
\label{div-rot-lem-alt}
\end{theorem}

Here, $\hmoom:=\hoc(\om)'$ and $\hmocom:=\hoom'$
and the distributional extensions 
\begin{align*}
\widetilde{\curl}:\ltom&\to\hmoom,
&
\widetilde{\div}:\ltom&\to\hmoom,\\
\widehat{\curl}:\ltom&\to\hmocom,
&
\widehat{\div}:\ltom&\to\hmocom
\end{align*}
of $\curl$ and $\div$, respectively, are defined for $E\in\ltom$ by
\begin{align*}
\widetilde{\curl}\,E\,(\Phi)&:=\scpltom{\curl\Phi}{E},
&
\Phi&\in\hoc(\om),\\
\widehat{\curl}\,E\,(\Phi)&:=\scpltom{\curl\Phi}{E},
&
\Phi&\in\hoom,\\
\widetilde{\div}\,E\,(\varphi)&:=-\scpltom{\na\varphi}{E},
&
\varphi&\in\hoc(\om),\\
\widehat{\div}\,E\,(\varphi)&:=-\scpltom{\na\varphi}{E},
&
\varphi&\in\hoom.
\end{align*}

Finally, we compare Theorem \ref{div-rot-lem-alt}
with the classical $\div$-$\curl$-lemma by Murat \cite{murat1978} and Tartar \cite{tartar1979},
which may be formulated as follows:

\begin{theorem}[classical $\div$-$\curl$-lemma]
Let $\om\subset\rt$ be an open set and let
\begin{itemize}
\item[\bf(i)]\quad
$E_{n},H_{n},E,H\in\ltom$, 
\item[\bf(i')]\quad
$E_{n}\wto E$ and $H_{n}\wto H$ in $\ltom$,
\item[\bf(ii)]\quad
$(\widetilde{\curl}E_{n})$ and $(\widetilde{\div}H_{n})$
be relatively compact in $\hmoom$.
\end{itemize}
Then
\begin{itemize}
\item[\bf(iii)]\quad
$\forall\,\varphi\in\cgen{}{\infty}{\ga}(\om)\qquad\scpltom{\varphi\,E_{n}}{H_{n}}\to\scpltom{\varphi\,E}{H}$.
\end{itemize}
\label{div-rot-lem-class}
\end{theorem}

\begin{acknowledgement}
The author is grateful to S\"oren Bartels for bringing up the topic of the $\div$-$\curl$-lemma,
and especially to Marcus Waurick for lots of inspiring discussions on the $\div$-$\curl$-lemma 
and for his substantial contributions to the Special Semester at RICAM in Linz late 2016.
\end{acknowledgement}

\bibliographystyle{spmpsci}
\bibliography{proc-wavesphen18-springer--pauly-bib}

\end{document}